\documentclass{article}
\usepackage{amssymb}
\usepackage{amsthm}
\usepackage{graphicx}
\usepackage{enumerate}
\usepackage{mathtools}
\usepackage{fullpage}
\usepackage{xcolor}

\usepackage{tikz}
\usetikzlibrary{patterns}
\usetikzlibrary{decorations.pathmorphing}
\usetikzlibrary{decorations.pathreplacing,decorations.markings}
\usetikzlibrary{shapes}

\makeatletter
\def\grd@save@target#1{%
	\def\grd@target{#1}}
\def\grd@save@start#1{%
	\def\grd@start{#1}}
\tikzset{
	grid with coordinates/.style={
		to path={%
			\pgfextra{%
				\edef\grd@@target{(\tikztotarget)}%
				\tikz@scan@one@point\grd@save@target\grd@@target\relax
				\edef\grd@@start{(\tikztostart)}%
				\tikz@scan@one@point\grd@save@start\grd@@start\relax
				\draw[minor help lines] (\tikztostart) grid (\tikztotarget);
				\draw[major help lines] (\tikztostart) grid (\tikztotarget);
				\grd@start
				\pgfmathsetmacro{\grd@xa}{\the\pgf@x/1cm}
				\pgfmathsetmacro{\grd@ya}{\the\pgf@y/1cm}
				\grd@target
				\pgfmathsetmacro{\grd@xb}{\the\pgf@x/1cm}
				\pgfmathsetmacro{\grd@yb}{\the\pgf@y/1cm}
				\pgfmathsetmacro{\grd@xc}{\grd@xa + \pgfkeysvalueof{/tikz/grid with coordinates/major step}}
				\pgfmathsetmacro{\grd@yc}{\grd@ya + \pgfkeysvalueof{/tikz/grid with coordinates/major step}}
				\foreach \x in {\grd@xa,\grd@xc,...,\grd@xb}
				\node[anchor=north] at (\x,\grd@ya) {\pgfmathprintnumber{\x}};
				\foreach \y in {\grd@ya,\grd@yc,...,\grd@yb}
				\node[anchor=east] at (\grd@xa,\y) {\pgfmathprintnumber{\y}};
			}
		}
	},
	minor help lines/.style={
		help lines,
		step=\pgfkeysvalueof{/tikz/grid with coordinates/minor step}
	},
	major help lines/.style={
		help lines,
		line width=\pgfkeysvalueof{/tikz/grid with coordinates/major line width},
		step=\pgfkeysvalueof{/tikz/grid with coordinates/major step}
	},
	grid with coordinates/.cd,
	minor step/.initial=.2,
	major step/.initial=1,
	major line width/.initial=0.25mm,
}

\pgfdeclarepatternformonly{north east lines wide}%
{\pgfqpoint{-1pt}{-1pt}}%
{\pgfqpoint{10pt}{10pt}}%
{\pgfqpoint{5pt}{5pt}}%
{
	\pgfsetlinewidth{.8pt}
	\pgfpathmoveto{\pgfqpoint{0pt}{0pt}}
	\pgfpathlineto{\pgfqpoint{9.1pt}{9.1pt}}
	\pgfusepath{stroke}
}

\tikzset{
	on each segment/.style={
		decorate,
		decoration={
			show path construction,
			moveto code={},
			lineto code={
				\path [#1]
				(\tikzinputsegmentfirst) -- (\tikzinputsegmentlast);
			},
			curveto code={
				\path [#1] (\tikzinputsegmentfirst)
				.. controls
				(\tikzinputsegmentsupporta) and (\tikzinputsegmentsupportb)
				..
				(\tikzinputsegmentlast);
			},
			closepath code={
				\path [#1]
				(\tikzinputsegmentfirst) -- (\tikzinputsegmentlast);
			},
		},
	},
	mid arrow/.style={postaction={decorate,decoration={
				markings,
				mark=at position .5 with {\arrow[#1]{stealth}}
	}}},
	rmid arrow/.style={postaction={decorate,decoration={
				markings,
				mark=at position .5 with {\arrowreversed[#1]{stealth}}
	}}},
	end arrow/.style={postaction={decorate,decoration={
				markings,
				mark=at position 1 with {\arrow[#1]{stealth}}
	}}},
	start arrow/.style={postaction={decorate,decoration={
				markings,
				mark=at position 0 with {\arrow[#1]{stealth}}
	}}},
	mid3 arrow/.style={postaction={decorate,decoration={
				markings,
				mark=at position .3 with {\arrow[#1]{stealth}}
	}}},
	rmid3 arrow/.style={postaction={decorate,decoration={
				markings,
				mark=at position .7 with {\arrowreversed[#1]{stealth}}
	}}},
	mid4 arrow/.style={postaction={decorate,decoration={
				markings,
				mark=at position .4 with {\arrow[#1]{stealth}}
	}}},
	rmid4 arrow/.style={postaction={decorate,decoration={
				markings,
				mark=at position .4 with {\arrowreversed[#1]{stealth}}
	}}},
}

\makeatother
\tikzset{every state/.style={minimum size=0pt}}

\tikzset{
	mark position/.style args={#1(#2)}{
		postaction={
			decorate,
			decoration={
				markings,
				mark=at position #1 with \coordinate (#2);
			}
		}
	}
}

\usetikzlibrary{3d}
\usepackage{xcolor}
\usepackage{xifthen}

\usetikzlibrary{decorations}

\pgfdeclaredecoration{ignore}{final}
{
	\state{final}{}
}

\pgfdeclaremetadecoration{middle}{initial}{
	\state{initial}[
	width={0pt},
	next state=middle
	]
	{\decoration{moveto}}
	
	\state{middle}[
	width={\pgfdecorationsegmentlength*\pgfmetadecoratedpathlength},
	next state=final
	]
	{\decoration{curveto}}
	
	\state{final}
	{\decoration{ignore}}
}

\tikzset{middle segment/.style={decoration={middle},decorate, segment length=#1}}

\newtheorem{thm}{Theorem}[section]
\newtheorem{lm}[thm]{Lemma}

\newtheorem{prop}[thm]{Proposition}

\newtheorem{cor}[thm]{Corollary}

\newtheorem{rmk}[thm]{Remark}
\numberwithin{equation}{section}

\newcommand{\bigO}{\mathcal{O}}

\newcommand{\complexC}{\mathbb{C}}

\newcommand{\ddbar}[2]{\frac{{\mathrm d}#1}{2\pi {\mathrm i}#2}}
\newcommand{\dg}{\mathrm{dist}}

\newcommand{\exe}{\mathrm{err}}

\newcommand{\ii}{\mathrm{i}}
\newcommand{\inn}{\mathrm{in}}

\newcommand{\limp}{\mathrm{p}}
\newcommand{\LL}{\mathrm{L}}
\newcommand{\out}{\mathrm{out}}
\newcommand{\prob}{\mathbb{P}}
\newcommand{\realR}{\mathbb{R}}

\newcommand{\RR}{\mathrm{R}}

\newcommand{\rz}{\mathrm{z}}

\begin{document}

\title{When the geodesic becomes rigid in the directed landscape}

\author{Zhipeng Liu\footnote{Department of Mathematics, University of Kansas, Lawrence, KS 66045. Email: \texttt{zhipeng@ku.edu}}}
\maketitle

\begin{abstract}
	When the value $L$ of the directed landscape at a point $(\mathbf{p};\mathbf{q})$ is sufficiently large, the geodesic from $\mathbf{p}$ to $\mathbf{q}$ is rigid and its location fluctuates of order $L^{-1/4}$  around its expectation. We further show that at a midpoint of the geodesic, the location of the geodesic and the value of the directed landscape after appropriate scaling converge to two independent Gaussians.
\end{abstract}

\section{Introduction}

Recently, there have been huge developments on the study of the Kardar-Parisi-Zhang universality class \cite{Baik-Deift-Johansson99,Johansson00,Johansson03,Borodin-Ferrari-Prahofer-Sasamoto07,Tracy-Widom08,Tracy-Widom09,Borodin-Corwin13,Matetski-Quastel-Remenik17,Dauvergne-Ortmann-Virag18,Johansson-Rahman19,Liu19}. Very recently, a four-parameter random field, the directed landscape, was constructed from one specific model in the class, the Brownian last passage percolation \cite{Dauvergne-Ortmann-Virag18}. It is believed that the directed landscape is the limiting law for all the models in the Kardar-Parisi-Zhang universality class, and this has been confirmed for several classic models \cite{Dauvergne-Virag21}. 

For the directed landscape, a lot of information is known, including the finite-dimensional distributions \cite{Matetski-Quastel-Remenik17,Johansson-Rahman19,Liu19}. On the other hand, less is known about the geodesic. There are studies on the properties of the geodesic very recently \cite{Basu-Sarkar-Sly17,Hammond20,Hammond-Sarkar20,Basu-Hoffman-Sly18,Basu-Ganguly-Hammon19,Bates-Ganguly-Hammond19b,Busani-Ferrari20,Dauvergne-Sarkar-Virag20,Corwin-Hammond-Hegde-Matetski21,Dauvergne-Virag21}. However, the explicit one-point distribution of the point-to-point geodesic was only obtained in \cite{Liu21}.

The goal of this paper is to investigate one property of the geodesic in the directed landscape using the formula obtained in \cite{Liu21}. 
We first introduce some notations and properties about the directed landscape and its geodesic. Let $\mathcal{L}(x,s;y,t)$ be the directed landscape. It satisfies the metric composition law
\begin{equation}
\mathcal{L}(x,r;y,t)=\max_{z\in\realR} \left(\mathcal{L}(x,r;z,s) + \mathcal{L}(z,s;y,t)\right)
\end{equation}
for any fixed $r<s<t$ and $x,y\in\realR$, and hence it has the reverse triangle inequality
\begin{equation}
\mathcal{L}(x,r;y,t)\ge \mathcal{L}(x,r;z,s) + \mathcal{L}(z,s;y,t).
\end{equation}
The geodesic in the directed landscape $\mathcal{L}$ from $(x,r)$ to $(y,t)$ is a continuous path $\pi=\pi_{x,r;y,t}$ which maximizes the following length with respect to $\mathcal{L}$
\begin{equation}
\label{eq:01}
\inf_{k\in\mathbb{N}} \inf_{r=s_0<s_1<\cdots<s_k=t} \sum_{i=1}^k\mathcal{L}(\pi(s_{i-1}),s_{i-1};\pi(s_i),s_i).
\end{equation}
It has been proved in \cite{Dauvergne-Ortmann-Virag18} that such a directed geodesic exists and is unique almost surely. Moreover, one can prove that the geodesic has the following properties
\begin{equation}
\label{eq:02}
\pi_{x,r;y,t}((1-s)r+st )\stackrel{law}=(t-r)^{2/3}\pi_{0,0;0,1}(s)+((1-s)x+sy),\quad s\in[0,1],
\end{equation}
where $\stackrel{law}=$ means that the two sides have the same law, and
\begin{equation}
\label{eq:03}
\begin{split}
&\mathcal{L}(\pi_{x,r;y,t}((1-s_1)r+s_1t),(1-s_1)r+s_1t;\pi_{x,r;y,t}((1-s_2)r+s_2t),(1-s_2)r+s_2t) \\
&\stackrel{law}= (t-r)^{1/3}\mathcal{L}(\pi_{0,0;0,1}(s_1),s_1;\pi_{0,0;0,1}(s_2),s_2)\\
&\quad-2(t-r)^{-1/3}(\pi_{0,0;0,1}(s_2)-\pi_{0,0;0,1}(s_1))(y-x)-\frac{s_2-s_1}{t-r}(y-x)^2
\end{split}
\end{equation}
for any $s_1, s_2$ satisfying $0\le s_1\le s_2\le 1$. In other words, the goedesic $\pi_{x,r;y,t}$ has the same law as $\pi_{0,0;0,1}$ after rescaling, and the directed landscape along the geodesic $\pi_{x,r;y,t}$ has the same law as along $\pi_{0,0;0,1}$ after shifting and rescaling. These two properties can be obtained straightforwardly using the formula~\eqref{eq:01} and the stationarity and rescaling properties of $\mathcal{L}$ described in \cite[Lemma 10.2]{Dauvergne-Ortmann-Virag18}. Hence we do not provide details here.

\bigskip

In this paper, we will fix the point $(0,0;0,1)$. Denote $\Pi(s)=\pi_{0,0;0,1}(s)$, $0\le s\le 1$, the geodesic from $(0,0)$ to $(0,1)$.  We also denote $\mathcal{L}(s)=\mathcal{L}(0,0;\Pi(s),s)$ for $0\le s\le 1$. We remark that the fact $\Pi(s)$ is on the geodesic implies $\mathcal{L}(s) + \mathcal{L}(\Pi(s),s;0,1)=\mathcal{L}(1)$.

The main result of this paper is about the fluctuations of $\Pi(s)$ and $\mathcal{L}(s)$ when $\mathcal{L}(1)=L$ becomes large. In the following theorem, the conditional probability $\prob\left(A\mid \mathcal{L}(1)=L\right)$ should be understood as $\lim_{\epsilon\to 0+}\prob\left(A\mid \mathcal{L}(1)-L\in(-\epsilon,\epsilon)\right)$, for an event $A$.

\begin{thm}[Rigidity of the geodesic]
	\label{thm:rigidity}
	For any fixed $s\in (0,1)$ and $x_1,x_2,\ell_1,\ell_2\in\realR$ satisfying $x_1<x_2$ and $\ell_1<\ell_2$, we have
	\begin{equation}
	\label{eq:thm1}
	\begin{split}
	&\lim_{L\to\infty}\prob\left(\frac{2{L}^{1/4}\Pi(s)}{\sqrt{s(1-s)}}\in (x_1,x_2), \frac{\mathcal{L}(s)-s{L}}{\sqrt{s(1-s)}{L}^{1/4}}\in (\ell_1,\ell_2) \:\Bigg\vert\: \mathcal{L}(1)= L\right)\\ &=\int_{x_1}^{x_2}p(x)\mathrm{d}x\int_{\ell_1}^{\ell_2}p(\ell)\mathrm{d}\ell,
	\end{split}
	\end{equation}
	where $p$ is the probability density function of the standard Gaussian distribution
	\begin{equation}
	p(x)=\frac{1}{\sqrt{2\pi}}e^{-\frac{x^2}{2}},\quad x\in\realR.
	\end{equation}
	Equivalently, the conditional joint density of $\frac{2\mathcal{L}^{1/4}\Pi(s)}{\sqrt{s(1-s)}}, \frac{\mathcal{L}(s)-s\mathcal{L}}{\sqrt{s(1-s)}\mathcal{L}^{1/4}}$ given $\mathcal{L}(1)=L$ converges to the product of two standard Gaussian densities as $L\to\infty$.
\end{thm}

\begin{rmk}
	Theorem~\ref{thm:rigidity} implies that $\Pi(s)$ becomes very rigid when $L$ becomes large. It fluctuates of order $L^{-1/4}$ and has Gaussian distribution after rescaling
	\begin{equation}
	\lim_{L\to\infty}\prob\left(\frac{2{L}^{1/4}\Pi(s)}{\sqrt{s(1-s)}}\in (x_1,x_2) \:\Bigg\vert\: \mathcal{L}(1)= L\right) =\int_{x_1}^{x_2}p(x)\mathrm{d}x.
	\end{equation}
	The rigidity of the geodesic conditioned on large $L$ is not surprising. For example in \cite{Basu-Ganguly19}, the authors were able to show that in the directed last passage percolation model with i.i.d. exponential weights, if the last passage time is $(4+\delta)n$, then the geodesic fluctuates of order $n^{1/2+o(1)}$ which is far smaller than the typical fluctuation order $n^{2/3}$ when the last passage time is about $4n+\bigO(n^{1/3})$. If we write $\delta = cn^{-2/3}L$, then $\Pi(s)\cdot n^{2/3}$ has an order $n^{1/2}$. So heuristically $\Pi(s)$ fluctuates as $\bigO(n^{-1/6})=\bigO(L^{-1/4})$ which is consistent with our result above. However, the limiting distribution of $\Pi(s)$ was not known to our best knowledge.
	
	Theorem~\ref{thm:rigidity} also implies that $\mathcal{L}(s)$ fluctuates of order $L^{1/4}$ and has Gaussian distribution after rescaling
	\begin{equation}
	\lim_{L\to\infty}\prob\left( \frac{\mathcal{L}(s)-s{L}}{{L}^{1/4}}\in (\ell_1,\ell_2) \:\Bigg\vert\: \mathcal{L}(1)= L\right) =\int_{\ell_1}^{\ell_2}p(\ell)\mathrm{d}\ell.
	\end{equation}
	Both the fluctuation order and the limiting distribution were not known to our best knowledge.
\end{rmk}

\begin{rmk}
	One can show that
	\begin{equation*}
	\begin{split}
	&\prob\left(\frac{2{L}^{1/4}\Pi(s)}{\sqrt{s(1-s)}}\in (x_1,x_2), \frac{\mathcal{L}(s)-s{L}}{\sqrt{s(1-s)}{L}^{1/4}}\in (\ell_1,\ell_2) \:\Bigg\vert\: \mathcal{L}(1)= L\right)\\ &=\int_{x_1}^{x_2}p(x)\mathrm{d}x\int_{\ell_1}^{\ell_2}p(\ell)\mathrm{d}\ell+L^{-3/4}\frac{2s-1}{2\sqrt{s(1-s)}}\int_{x_1}^{x_2}p(x)\mathrm{d}x\int_{\ell_1}^{\ell_2}\ell p(\ell)\mathrm{d}\ell+\bigO(L^{-3/2})
	\end{split}
	\end{equation*}
	when $L$ becomes large. The error terms can be explicitly evaluated by a more careful calculation using the same argument of this paper. However, the calculation is quite cumbersome so we do not include it here.
\end{rmk}

By combing the equations~\eqref{eq:02} and~\eqref{eq:03}, and Theorem~\ref{thm:rigidity}, we have the following general result for the geodesic from $(x,r)$ to $(y,t)$, conditioned on $\mathcal{L}(x,r;y,t)$ goes to infinity.
\begin{cor}
	Suppose $x,y\in\realR$, and $r,t\in\realR_{\ge 0}$ are all fixed. Assume that $r<t$. Conditioned on $\mathcal{L}(x,r;y,t)=L$ goes to infinity, we have
	\begin{equation*}
	\frac{2L^{1/4}\left(\pi_{x,r;y,t}((1-s)r+st)- ((1-s)x+sy)\right)}
	{(t-r)^{3/4}\sqrt{s(1-s)}}
	\end{equation*}
	and
	\begin{equation*}
	\frac{\mathcal{L}\left(x,r;\pi_{x,t;y,t}((1-s)r+st),(1-s)r+st\right)-sL}{(t-r)^{1/4}\sqrt{s(1-s)}L^{1/4}}
	\end{equation*}
	converge to two independent standard Gaussians in distribution.
\end{cor}
\bigskip

Our approach relies on the explicit joint density function of $\mathcal{L}(s)=\mathcal{L}(0,0;\Pi(s),s), \mathcal{L}(1)-\mathcal{L}(s)=\mathcal{L}(\Pi(s),s;0,1)$ and $\Pi(s)$ which was obtained very recently \cite{Liu21}. As we mentioned at the beginning of this paper, \cite{Liu21} gives the first explicit formula of the one-point distribution of the point-to-point geodesic. However, the formula in \cite{Liu21} looks very complicated and people might doubt whether such a formula has any probabilistic applications. This paper provides a simple but beautiful application.

\bigskip

The organization of the paper is as follows. In Section~\ref{sec:joint_density}, we introduce the joint density function mentioned above. Then we prove Theorem~\ref{thm:rigidity} In Section~\ref{sec:proof}.

\section*{Acknowledgments}
We would like to thank Jinho Baik, Ivan Corwin, Duncan Dauvergne, and B\'alint Vir\'ag for the comments and suggestions. The work was supported by the University of Kansas Start Up Grant, the University of Kansas New Faculty General Research Fund, Simons Collaboration Grant No. 637861, and  NSF grant
DMS-1953687. 
\section{The joint density function}
\label{sec:joint_density}

In this section, we introduce the joint density function of $\mathcal{L}(s), \mathcal{L}(1)-\mathcal{L}(s)$ and $\Pi(s)$, which was obtained in \cite[Corollary 1.6]{Liu21}. We rephrase it as follows.
\begin{prop}[\cite{Liu21}]
	\label{prop:joint_density}
	The joint density function of $\mathcal{L}(s), \mathcal{L}(1)-\mathcal{L}(s)$ and $\Pi(s)$ is
	\begin{equation}
	2\limp\left(\ell_1+\frac{x^2}{s},\ell_2+\frac{x^2}{1-s},2x;s\right),
	\end{equation}
	where $\limp(\ell_1,\ell_2,x;s)$ is the joint density function introduced in \cite{Liu21}, and as described in~\eqref{eq:def_limiting_density} below.
\end{prop}

\begin{figure}[t]
	\centering
	\begin{tikzpicture}[scale=1]
	\draw [line width=0.4mm,lightgray] (-3,0)--(3,0) node [pos=1,right,black] {$\realR$};
	\draw [line width=0.4mm,lightgray] (0,-1)--(0,1) node [pos=1,above,black] {$\mathrm{i}\realR$};
	\fill (0,0) circle[radius=2pt] node [below,shift={(0pt,-3pt)}] {$0$};

	\path [draw=black,thick,postaction={mid3 arrow={black,scale=1.2}}]	(3,-5.1/6) 
	to [out=160,in=-90] (2,0)
	to [out=90,in=-160] (3,5.1/6);
	\path [draw=black,thick,postaction={mid3 arrow={black,scale=1.2}}]	(2.6,-6.1/6) 
	to [out=160,in=-90] (1.5,0)
	to [out=90,in=-160] (2.6,6.1/6);
	\path [draw=black,thick,postaction={mid3 arrow={black,scale=1.2}}]	(2.3,-7.1/6) 
	to [out=160,in=-90] (1,0)
	to [out=90,in=-160] (2.3,7.1/6);
	
	\path [draw=black,thick,postaction={mid3 arrow={black,scale=1.5}}]	(-3,-5.1/6) 
	to [out=20,in=-90] (-2,0)
	to [out=90,in=-20] (-3,5.1/6);
	\path [draw=black,thick,postaction={mid3 arrow={black,scale=1.5}}]	(-2.6,-6.1/6) 
	to [out=20,in=-90] (-1.5,0)
	to [out=90,in=-20] (-2.6,6.1/6);
	\path [draw=black,thick,postaction={mid3 arrow={black,scale=1.5}}]	(-2.3,-7.1/6) 
	to [out=20,in=-90] (-1,0)
	to [out=90,in=-20] (-2.3,7.1/6);
	\end{tikzpicture}
	\caption{The three contours in the left half plane from left to right are $\Gamma_{\LL,\inn}, \Gamma_\LL$ and $\Gamma_{\LL,\out}$ respectively, and the three contours in the right half plane from left to right are $\Gamma_{\RR,\out}, \Gamma_{\RR}$ and $\Gamma_{\RR,\inn}$ respectively.}\label{fig:contours_limit}
\end{figure}
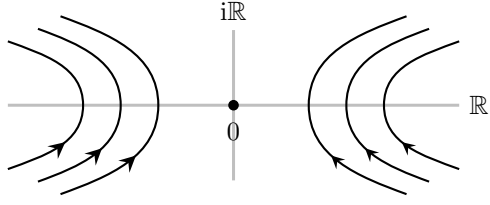

We need to introduce six contours before introducing the function $\limp$. Suppose $\Gamma_{\LL,\inn},\Gamma_{\LL}$ and $\Gamma_{\LL,\out}$ are three disjoint contours on the left half plane each of which starts from $e^{-2\pi\ii/3}\infty$ and ends to $e^{2\pi\ii/3}\infty$. Here $\Gamma_{\LL,\inn}$ is the leftmost contour and $\Gamma_{\LL,\out}$ is the rightmost contour. The index ``$\inn$'' and ``$\out$'' refer to the relative location compared with $-\infty$. Similarly, suppose  $\Gamma_{\RR,\inn},\Gamma_{\RR}$ and $\Gamma_{\RR,\out}$ are three disjoint contours on the right half plane each of which starts from $e^{-\pi\ii/3}\infty$ and ends to $e^{\pi\ii/3}\infty$. Here the index ``$\inn$'' and ``$\out$'' refer to the relative location compared with $+\infty$, hence $\Gamma_{\RR,\inn}$ is the rightmost contour and $\Gamma_{\RR,\out}$ is the leftmost contour. See Figure~\ref{fig:contours_limit} for an illustration of these contours.

\bigskip

The probability density function $\limp(\ell_1,\ell_2,x;s)$ is defined to be
\begin{equation}
\label{eq:def_limiting_density}
\limp(\ell_1,\ell_2,x;s)
:=\oint_0\ddbar{\mathrm{z}}{\mathrm{(1-z)^2}}\sum_{k_1,k_2\ge 1}
\frac{1}{(k_1!k_2!)^2}
\mathrm{T}_{k_1,k_2}(\mathrm{z};\ell_1,\ell_2,x;s)
\end{equation}
with
\begin{equation}
\label{eq:def_limit_T}
\begin{split}
&\mathrm{T}_{k_1,k_2}(\mathrm{z};\ell_1,\ell_2,x;s)\\
&:=
\prod_{i_1=1}^{k_1}     
\left(  \frac{1}{1-\mathrm{z}} \int_{\Gamma_{\LL,\inn}} 	    
\ddbar{\xi^{(1)}_{i_1}}{}   
-\frac{\mathrm{z}}{1-\mathrm{z}} \int_{\Gamma_{\LL,\out}} \ddbar{\xi^{(1)}_{i_1}}{}
\right)
\left(  \frac{1}{1-\mathrm{z}}  \int_{\Gamma_{\RR,\inn}}  
\ddbar{\eta^{(1)}_{i_1}}{}
-\frac{\mathrm{z}}{1-\mathrm{z}}   \int_{\Gamma_{\RR,\out}}
\ddbar{\eta^{(1)}_{i_1}}{}
\right)
\\
&\quad \cdot 
\prod_{i_2=1}^{k_2}  
\int_{\Gamma_\LL} \ddbar{\xi^{(2)}_{i_2}}{} 
\int_{\Gamma_\RR} \ddbar{\eta^{(2)}_{i_2}}{} 
\cdot    
\left(1-\mathrm{z}\right)^{k_2}
\left(1-\frac{1}{\mathrm{z}}\right)^{k_1}
\cdot 
\frac
{\mathrm{f}_1(\boldsymbol{\xi}^{(1)};\ell_1)
	\mathrm{f}_2(\boldsymbol{\xi}^{(2)};\ell_2)}
{\mathrm{f}_1(\boldsymbol{\eta}^{(1)};\ell_1)
	\mathrm{f}_2(\boldsymbol{\eta}^{(2)};\ell_2)}
\cdot
\mathrm{H}(\boldsymbol{\xi}^{(1)},\boldsymbol{\eta}^{(1)};
\boldsymbol{\xi}^{(2)},\boldsymbol{\eta}^{(2)})\\
&\quad \cdot 
\prod_{\ell=1}^2
\frac{\left( \Delta(\boldsymbol{\xi}^{(\ell)}) \right)^2
	\left( \Delta(\boldsymbol{\eta}^{(\ell)}) \right)^2
}
{\left( \Delta(\boldsymbol{\xi}^{(\ell)}; \boldsymbol{\eta}^{(\ell)}) \right)^2}
\cdot 
\frac{\Delta(\boldsymbol{\xi}^{(1)};\boldsymbol{\eta}^{(2)})
	\Delta(\boldsymbol{\eta}^{(1)};\boldsymbol{\xi}^{(2)})
}
{\Delta(\boldsymbol{\xi}^{(1)};\boldsymbol{\xi}^{(2)})
	\Delta(\boldsymbol{\eta}^{(1)};\boldsymbol{\eta}^{(2)})
}
\end{split}
\end{equation}
where the vectors $\boldsymbol{\xi}^{(\ell)}=(\xi_1^{(\ell)},\cdots,\xi_{i_\ell}^{(\ell)})$ and $\boldsymbol{\eta}^{(\ell)}=(\eta_1^{(\ell)},\cdots,\eta^{(\ell)}_{i_\ell})$ for $\ell\in\{1,2\}$, the functions $\mathrm{f}_1,\mathrm{f}_2$ are defined by
\begin{equation}
\label{eq:def_rmf}
\begin{split}
\mathrm{f}_1(\zeta;\ell_1)&:=\exp\left(-\frac{s}{3}\zeta^3-\frac12\mathrm{x}\zeta^2+\left(\ell_1-\frac{{x}^2}{4s}\right)\zeta\right),\\
\mathrm{f}_2(\zeta;\ell_2)& := \exp\left(-\frac{(1-s)}{3}\zeta^3+\frac12\mathrm{x}\zeta^2+\left(\ell_2-\frac{{x}^2}{4(1-s)}\right)\zeta\right),
\end{split}
\end{equation}
and the function $\mathrm{H}$ is defined by
\begin{equation}
\label{eq:def_rmH}
\mathrm{H}(\boldsymbol{\xi}^{(1)},\boldsymbol{\eta}^{(1)};
\boldsymbol{\xi}^{(2)},\boldsymbol{\eta}^{(2)})
=\frac{1}{12}\mathrm{S}_1^4+\frac14\mathrm{S}_2^2-\frac13\mathrm{S}_1\mathrm{S}_3
\end{equation}
with
\begin{equation}
\label{eq:def_S}
\mathrm{S}_j
=\mathrm{S}_j(\boldsymbol{\xi}^{(1)},\boldsymbol{\eta}^{(1)};
\boldsymbol{\xi}^{(2)},\boldsymbol{\eta}^{(2)})
=\sum_{i_1=1}^{k_1}
\left( \left(\xi_{i_1}^{(1)}\right)^{j}
-\left(\eta_{i_1}^{(1)}\right)^{j}
\right)
-\sum_{i_2=1}^{k_2}
\left( \left(\xi_{i_2}^{(2)}\right)^{j}
-\left(\eta_{i_2}^{(2)}\right)^{j}
\right).
\end{equation}
We also used the notations
\begin{equation*}
\Delta(W):=\prod_{1\le i<j\le n}(w_j-w_i), \quad \Delta(W;W'):=\prod_{i=1}^n\prod_{i'=1}^{n'}(w_i-w'_{i'}),\quad f(W)=\prod_{i=1}^n f(w)
\end{equation*}
for any two vectors $W=(w_1,\cdots,w_n)\in\complexC^n$ and $W'=(w'_1,\cdots,w'_{n'})\in\complexC^{n'}$ and any function $f:\complexC\to\complexC$. Here we allow the empty product and set it to be $1$.

\section{Proof of Theorem~\ref{thm:rigidity}}
\label{sec:proof}
It is well known that $\mathcal{L}(1)$ has the GUE Tracy-Widom distribution $F_{GUE}(L)$, which can be defined by
\begin{equation}
F_{GUE}(L)=\exp\left(-\int_L^\infty \left(\ell-L\right)u^2(\ell)\mathrm{d}\ell\right)
\end{equation}
in terms of the  Hastings-McLeod  solution $u=u(\ell)$ to the Painlev\'e-II equation 
\begin{equation}
u''= \ell u + 2u^3,\qquad u(\ell)\to \frac{1}{2\sqrt{\pi}\ell^{1/4}}e^{-\frac{2}{3}\ell^{3/2}}(1+o(1)),\quad \text{as}\quad \ell\to\infty.
\end{equation}
It is known from the above formulas (see, for example \cite[equation (25)]{Baik-Buckingham-DiFranco08} for the expansion of the squared root of ${F_{GUE}(L)}$) that $F_{GUE}(L)$ has the following right tail behavior
\begin{equation}
F_{GUE}(L) \approx 1-\frac{1}{16\pi L^{3/2}}e^{-\frac{4}{3}L^{3/2}},\quad \text{as  } L\to\infty,
\end{equation}
and its density function $f_{GUE}(L)$ has
\begin{equation}
\label{eq:aux_02}
f_{GUE}(L) \approx\frac{1}{8\pi L}e^{-\frac{4}{3}L^{3/2}}, \quad \text{as  } L\to\infty.
\end{equation}
Here and afterwards, we use $\approx$ to denote the leading order asymptotics if the next order term is not needed in our analysis.

Now we come to the proof of Theorem~\ref{thm:rigidity}. For fixed $s\in(0,1)$ and $x,\ell\in\realR$, we denote
\begin{equation}
\label{eq:scaling}
\begin{split}
X=X(L)&:=\frac{x\sqrt{s(1-s)}}{L^{1/4}},\\
L_1=L_1(L)&:=sL+\sqrt{s(1-s)}L^{1/4}\ell+\frac{x^2(1-s)}{4L^{1/2}},\\ L_2=L_2(L)&:=(1-s)L-\sqrt{s(1-s)}L^{1/4}\ell+\frac{x^2s}{4L^{1/2}}.
\end{split}
\end{equation}
In order to show Theorem~\ref{thm:rigidity}, it is sufficient to prove, by using Proposition~\ref{prop:joint_density}, the following asymptotic result holds uniformly for $x$ and $\ell$ on compact intervals
\begin{equation}
\label{eq:aux_01}
\lim_{L\to\infty}\frac{s(1-s)\limp\left(L_1,L_2,X;s\right)}{f_{GUE}(L)}=\frac{1}{\sqrt{2\pi}}e^{-\frac{1}{2}\ell^2}\cdot\frac{1}{\sqrt{2\pi}} e^{-\frac{1}{2}x^2}
\end{equation}
where $f_{GUE}(L)$ is the probability density function of $\mathcal{L}$.

\bigskip

Inserting~\eqref{eq:aux_02} in~\eqref{eq:aux_01}, we find that~\eqref{eq:aux_01} is equivalent to
\begin{equation}
\label{eq:aux_03}
\limp(L_1,L_2,X;s)=\frac{1}{16\pi^2 s(1-s) L}e^{-\frac{4}{3}L^{3/2}}\left(e^{-\frac{1}{2}(\ell^2+x^2)}+o(1)\right),\quad \text{as}\quad L\to\infty
\end{equation}
uniformly for $x$ and $\ell$ on compact intervals.

Recall the definition~\eqref{eq:def_limiting_density} of $\limp$. We have the following two lemmas.
\begin{lm}
	\label{lm:01}
	With the scaling~\eqref{eq:scaling}, we have
	\begin{equation}
	\oint_0\ddbar{\rz}{(1-\rz)^2}\mathrm{T}_{1,1}(\rz;L_1,L_2,X;s)=\frac{1}{16\pi^2 s(1-s) L}e^{-\frac{4}{3}L^{3/2}}\left(e^{-\frac{1}{2}(\ell^2+x^2)}+o(1)\right)
	\end{equation}
	as $L\to\infty$, where the error term $o(1)$ holds uniformly for $x$ and $\ell$ on compact intervals.
\end{lm}
\begin{lm}
	\label{lm:02}
	Suppose $|z|=r\in (0,1)$ is fixed. With the scaling~\eqref{eq:scaling}, for any fixed $\epsilon\in(0,1)$, there exists a constant $C>0$ such that for sufficiently large $L$, we have
	\begin{equation}
	\left|\mathrm{T}_{k_1,k_2}(\rz;L_1,L_2,X;s)\right|\le k_1^{k_1/2}k_2^{k_2/2}(k_1+k_2)^{(k_1+k_2)/2}C^{k_1+k_2} e^{-\frac{4(1-\epsilon)}{3}(sk_1+(1-s)k_2)L^{3/2}}
	\end{equation}
	for all $k_1,k_2\ge 1$.
\end{lm}

The proof of these two lemmas will be given in the following two subsections. Now we use these two lemmas to show~\eqref{eq:aux_03}. Using~\eqref{eq:def_limiting_density}, we write
\begin{equation}
\label{eq:aux_04}
\begin{split}
\limp(L_1,L_2,X;s)&=\oint_0\ddbar{\rz}{(1-\rz)^2}\mathrm{T}_{1,1}(\rz;L_1,L_2,X;s)\\
&\quad+\oint_{|\rz|=r}\ddbar{\rz}{(1-\rz)^2}\sum_{\substack{k_1+k_2\ge 3\\ k_1,k_2\ge 1}}\frac{1}{(k_1!k_2!)^2}\mathrm{T}_{k_1,k_2}(\rz;L_1,L_2,X;s).
\end{split}
\end{equation}
By Lemma~\ref{lm:01}, the first term on the right hand side of~\eqref{eq:aux_04} gives the desired leading term in~\eqref{eq:aux_03}. We then apply Lemma~\ref{lm:02} for the second term on the right hand side of~\eqref{eq:aux_04} and obtain
\begin{equation}
\begin{split}
&\left|\oint_{|\rz|=r}\ddbar{\rz}{(1-\rz)^2}\sum_{\substack{k_1+k_2\ge 3\\ k_1,k_2\ge 1}}\frac{1}{(k_1!k_2!)^2}\mathrm{T}_{k_1,k_2}(\rz;L_1,L_2,X;s)\right|\\
&\le \oint_{|\rz|=r}\frac{|\mathrm{d}\rz|}{2\pi|1-\rz|^2}\sum_{\substack{k_1+k_2\ge 3\\ k_1,k_2\ge 1}}\frac{1}{(k_1!k_2!)^2}k_1^{k_1/2}k_2^{k_2/2}(k_1+k_2)^{(k_1+k_2)/2}C^{k_1+k_2} e^{-\frac{4(1-\epsilon)}{3}(sk_1+(1-s)k_2)L^{3/2}}\\
&\le C'\max_{\substack{ k_1+k_2\ge 3\\ k_1,k_2\ge 1}}e^{-\frac{4(1-\epsilon)}{3}(sk_1+(1-s)k_2)L^{3/2}}
\end{split}
\end{equation}
for sufficiently large $L$, where $C'$ is a constant given by
\begin{equation*}
C'=\oint_{|\rz|=r}\frac{|\mathrm{d}\rz|}{2\pi|1-\rz|^2}\sum_{\substack{k_1+k_2\ge 3\\ k_1,k_2\ge 1}}\frac{1}{(k_1!k_2!)^2}k_1^{k_1/2}k_2^{k_2/2}(k_1+k_2)^{(k_1+k_2)/2}C^{k_1+k_2}.
\end{equation*}
Note that the above sum in the definition of $C'$ is absolutely convergent using the Stirling's approximation formula. On the other hand, we can choose $\epsilon$ small enough such that $\frac{4(1-\epsilon)}{3}(sk_1+(1-s)k_2)>\frac{4}{3}$ for all $k_1+k_2\ge 3$ and $k_1,k_2\ge 1$. Then the second term on the right hand side of~\eqref{eq:aux_04} is negligible compared to the first term as $L$ becomes large. Hence we complete the proof of~\eqref{eq:aux_03}.

\bigskip

It remains to prove the two lemmas. We provide the proof in the following two subsections.

\subsection{Proof of Lemma~\ref{lm:01}}
In the proof below, we fixed $x$ and $\ell$ for simplification. However, it is easy to check that the proof is also valid uniformly for $x$ and $\ell$ on compact intervals since our error terms depend on $x$ and $\ell$ continuously.

Since $k_1=k_2=1$, we drop the subscript $1$ in the integration variables for simplification
\begin{equation}
\xi^{(1)}=\xi_1^{(1)},\quad \eta^{(1)}= \eta_1^{(1)},\quad \xi^{(2)}=\xi_1^{(2)}, \quad\text{and}\quad \eta^{(2)}=\eta_1^{(2)}.
\end{equation}
Recall the formula for $\mathrm{T}_{1,1}$ in~\eqref{eq:def_limit_T}. We realize that if we expand the integrals, there are four possible combinations in $\mathrm{T}_{1,1}$
\begin{equation}
\begin{split}
&-\frac{1}{\rz}\int_{\Gamma_{\LL,\inn}}\int_{\Gamma_{\RR,\inn}}\int_{\Gamma_\LL}\int_{\Gamma_\RR},
\qquad\frac{1}{1-\rz}\int_{\Gamma_{\LL,\out}}\int_{\Gamma_{\RR,\inn}}\int_{\Gamma_\LL}\int_{\Gamma_\RR},\\
&\frac{1}{1-\rz}\int_{\Gamma_{\LL,\inn}}\int_{\Gamma_{\RR,\out}}\int_{\Gamma_\LL}\int_{\Gamma_\RR},
\qquad-\rz\int_{\Gamma_{\LL,\out}}\int_{\Gamma_{\RR,\out}}\int_{\Gamma_\LL}\int_{\Gamma_\RR},
\end{split}
\end{equation}
where we ignore the integrand for simplification. Now if we take the $\rz$ integration $\oint_0\ddbar{\rz}{(1-\rz)^2}$ for these four combinations, clearly only the first one survives
\begin{equation}
\label{eq:aux_05}
\oint_0\ddbar{\rz}{(1-\rz)^2}\left(-\frac{1}{\rz}\int_{\Gamma_{\LL,\inn}}\int_{\Gamma_{\RR,\inn}}\int_{\Gamma_\LL}\int_{\Gamma_\RR}\right)=-\int_{\Gamma_{\LL,\inn}}\int_{\Gamma_{\RR,\inn}}\int_{\Gamma_\LL}\int_{\Gamma_\RR}.
\end{equation}
Now we write down the explicit formula for the above integrals, also note
\begin{equation}
\mathrm{H}(\boldsymbol{\xi}^{(1)},\boldsymbol{\eta}^{(1)};\boldsymbol{\xi}^{(2)},\boldsymbol{\eta}^{(2)})=(\xi^{(1)}-\eta^{(1)})(\xi^{(2)}-\eta^{(2)})(\eta^{(1)}-\eta^{(2)})(\xi^{(1)}-\xi^{(2)})
\end{equation}
which follows from a direct calculation. We obtain
\begin{equation}
\label{eq:aux_06}
\begin{split}
&\oint_0\ddbar{\rz}{(1-\rz)^2}\mathrm{T}_{1,1}(\rz;L_1,L_2,X;s)\\
&=-\int_{\Gamma_{\LL,\inn}}\ddbar{\xi^{(1)}}{}\int_{\Gamma_{\RR,\inn}}\ddbar{\eta^{(1)}}{}\int_{\Gamma_\LL}\ddbar{\xi^{(2)}}{}\int_{\Gamma_\RR}\ddbar{\eta^{(2)}}{}\frac{e^{g_1(\xi^{(1)})}e^{g_2(\xi^{(2)})}}{e^{g_1(\eta^{(1)})}e^{g_2(\eta^{(2)})}}\frac{\left(\xi^{(1)}-\eta^{(2)}\right)\left(\eta^{(1)}-\xi^{(2)}\right)}{\left(\xi^{(1)}-\eta^{(1)}\right)\left(\xi^{(2)}-\eta^{(2)}\right)}\\
\end{split}
\end{equation}
where
\begin{equation}
\label{eq:def_g}
g_1(\xi)=-\frac{s}{3}\xi^3-\frac{1}{2}X\xi^2+\left(L_1-\frac{X^2}{4s}\right)\xi,\quad g_2(\xi)=-\frac{1-s}{3}\xi^3+\frac{1}{2}X\xi^2+\left(L_2-\frac{X^2}{4(1-s)}\right)\xi.
\end{equation}

Now we apply a standard steepest descent analysis to~\eqref{eq:aux_06}. Note that if we only consider the leading term, $g_1(\xi)\approx s\left(-\xi^3/3 + L\xi\right)$ and $g_2(\xi)\approx (1-s)(-\xi^3/3+L\xi)$ by~\eqref{eq:scaling}. Thus the main contribution of the integral comes from $\xi^{(1)},\xi^{(2)}\approx -\sqrt{L}$ and $\eta^{(1)},\eta^{(2)}\approx \sqrt{L}$. To further analyze the integral we need to zoom in a neighborhood of these two points.

We deform the contours $\Gamma_{*}$, with $*=\{\LL,\inn\},\{\LL\},\{\RR,\inn\}$, and $\{\RR\}$, in the following way. Here the braces are just introduced to avoid notation confusion but they should not appear in the subscripts. For example, $\Gamma_*=\Gamma_{\LL,\inn}$ if $*=\{\LL,\inn\}$. 

The contour $\Gamma_*$ has a vertical part $\Gamma_{*}^{(L)}$
\begin{equation}
\begin{split}
&\Gamma_{*}^{(L)}\\
&:=\begin{dcases}
-\frac{X}{2s}-\sqrt{L}-\frac{1}{2}\sqrt{\frac{1-s}{s}}\ell L^{-1/4} +L^{-1/4}\cdot \left\{\mathrm{i}u^{(1)}: -\log L\le u^{(1)}\le \log L \right\}, & *=\{\LL,\inn\},\\
\frac{X}{2(1-s)}-\sqrt{L} +\frac{1}{2}\sqrt{\frac{s}{1-s}}\ell L^{-1/4}+L^{-1/4}\cdot \left\{\mathrm{i}u^{(2)}: -\log L\le u^{(2)}\le \log L \right\}, & *=\{\LL\},\\
-\frac{X}{2s}+\sqrt{L} +\frac{1}{2}\sqrt{\frac{1-s}{s}}\ell L^{-1/4}+L^{-1/4}\cdot \left\{\mathrm{i}v^{(1)}: -\log L\le v^{(1)}\le \log L \right\}, & *=\{\RR,\inn\},\\
\frac{X}{2(1-s)}+\sqrt{L} -\frac{1}{2}\sqrt{\frac{s}{1-s}}\ell L^{-1/4}+L^{-1/4}\cdot \left\{\mathrm{i}v^{(2)}: -\log L\le v^{(2)}\le \log L \right\}, & *=\{\RR\},
\end{dcases} 
\end{split}
\end{equation}
and a second part $\Gamma_{*}^{(\exe)}$ containing two rays
\begin{equation}
\Gamma_{*}^{(\exe)}:=\begin{dcases}
-\frac{X}{2s}-\sqrt{L}-\frac{1}{2}\sqrt{\frac{1-s}{s}}\ell L^{-1/4}\pm \mathrm{i}L^{-1/4}\log L+ e^{\pm 2\pi\mathrm{i}/3}\cdot\realR_{\ge 0},&*=\{\LL,\inn\},\\
\frac{X}{2(1-s)}-\sqrt{L}+\frac{1}{2}\sqrt{\frac{s}{1-s}}\ell L^{-1/4}\pm \mathrm{i}L^{-1/4}\log L+ e^{\pm 2\pi\mathrm{i}/3}\cdot\realR_{\ge 0},&*=\{\LL\},\\
-\frac{X}{2s}+\sqrt{L}+\frac{1}{2}\sqrt{\frac{1-s}{s}}\ell L^{-1/4}\pm \mathrm{i}L^{-1/4}\log L+ e^{\pm \pi\mathrm{i}/3}\cdot\realR_{\ge 0},&*=\{\RR,\inn\},\\
\frac{X}{2(1-s)}+\sqrt{L}-\frac{1}{2}\sqrt{\frac{s}{1-s}}\ell L^{-1/4}\pm \mathrm{i}L^{-1/4}\log L+ e^{\pm \pi\mathrm{i}/3}\cdot\realR_{\ge 0},&*=\{\RR\}.
\end{dcases}
\end{equation}
Note that the contours $\Gamma_{\LL,\inn}$ and $\Gamma_{\RR,\inn}$ do not intersect during the deformation, similarly $\Gamma_{\LL}$ and $\Gamma_{\RR}$ do not intersect. Thus the denominator in the integrand of~\eqref{eq:aux_06} does not vanish during the deformation of these contours.

Now we estimate $g_i$ on $\Gamma_*^{(L)}$ and $\Gamma_*^{(\exe)}$ respectively.

For 
\begin{equation}
\zeta=-\frac{X}{2s}\mp\sqrt{L}\mp\frac{1}{2}\sqrt{\frac{1-s}{s}}\ell L^{-1/4}+ L^{-1/4}\mathrm{i}w,\quad-\log L\le w\le \log L,
\end{equation} we have, after inserting~\eqref{eq:scaling} and using Taylor expansion,
\begin{equation}
\label{eq:aux_07}
\begin{split}
&\pm\left(g_1(\zeta)-\left(\frac{X^3}{24s^2}-\frac{XL_1}{2s}\right)\right)\\
&=\mp\frac{s}{3}\left(\zeta+\frac{X}{2s}\right)^3\pm L_1\left(\zeta+\frac{X}{2s}\right)\\
&=-\frac{2}{3}sL^{3/2}-\sqrt{s(1-s)}\ell L^{3/4} -\frac{1}{4}(1-s)\ell^2-\frac{1}{4}(1-s)x^2 -sw^2+\bigO(L^{-3/4}(\log L)^3).
\end{split}
\end{equation}
For 
\begin{equation}
\zeta=-\frac{X}{2s}\mp\sqrt{L}\mp\frac{1}{2}\sqrt{\frac{1-s}{s}}\ell L^{-1/4}+ L^{-1/4}\mathrm{i}\log L +(\mp1+ \sqrt{3}\mathrm{i})w,\quad w\in\realR_{\ge 0}, 
\end{equation}
or
\begin{equation}
\zeta=-\frac{X}{2s}\mp\sqrt{L}\mp\frac{1}{2}\sqrt{\frac{1-s}{s}}\ell L^{-1/4}- L^{-1/4}\mathrm{i}\log L +(\mp1- \sqrt{3}\mathrm{i})w,\quad w\in\realR_{\ge 0}, 
\end{equation}
we have, after a tedious but straightforward calculation,
\begin{equation}
\label{eq:aux_08}
\begin{split}
&\pm\mathrm{Re}\left(g_1(\zeta)-\left(\frac{X^3}{24s^2}-\frac{XL_1}{2s}\right)\right)\\
&=\mathrm{Re}\left(\mp\frac{s}{3}\left(\zeta+\frac{X}{2s}\right)^3\pm L_1\left(\zeta+\frac{X}{2s}\right)\right)\\
&=-\frac{2}{3}sL^{3/2}-\sqrt{s(1-s)}\ell L^{3/4} -\frac{1}{4}(1-s)\ell^2-\frac{1}{4}(1-s)x^2 -s(\log L)^2+\bigO(L^{-3/4}(\log L)^3)\\
&\quad -\frac{8s}{3}w^3 -(2sL^{1/2}+\sqrt{s(1-s)}\ell L^{-1/4} +2\sqrt{3}sL^{-1/4}\log L)w^2\\
&\quad
-\left(2\sqrt{3}sL^{1/4}\log L+sL^{-1/2}(\log L)^2 +\sqrt{3s(1-s)}\ell L^{-1/2}\log L\right.\\
&\qquad\left. +\frac{x^2(1-s)}{4}L^{-1/2}-\frac{1-s}{4}\ell^2 L^{-1/2}\right)w,
\end{split}
\end{equation}
where the error term is independent of $w$. Note that the $w$ terms behaves like $-\frac{8s}{3}w^3 -2sL^{1/2}w^2 -2\sqrt{3}sL^{1/4}\log L w$ as $L$ becomes large.

Similarly, for
\begin{equation}
\zeta=\frac{X}{2(1-s)}\mp\sqrt{L}\pm\frac{1}{2}\sqrt{\frac{s}{1-s}}\ell L^{-1/4}+ L^{-1/4}\mathrm{i}w,\quad-\log L\le w\le \log L,
\end{equation} we have
\begin{equation}
\label{eq:aux_09}
\begin{split}
&\pm\left(g_2(\zeta)+\left(\frac{X^3}{24(1-s)^2}-\frac{XL_2}{2(1-s)}\right)\right)\\
&=\mp\frac{1-s}{3}\left(\zeta-\frac{X}{2(1-s)}\right)^3\pm L_2\left(\zeta-\frac{X}{2(1-s)}\right)\\
&=-\frac{2}{3}(1-s)L^{3/2}+\sqrt{s(1-s)}\ell L^{3/4} -\frac{1}{4}s\ell^2-\frac{1}{4}sx^2 -(1-s)w^2+\bigO(L^{-3/4}(\log L)^3).
\end{split}
\end{equation}
For 
\begin{equation}
\zeta=\frac{X}{2(1-s)}\mp\sqrt{L}\pm\frac{1}{2}\sqrt{\frac{1-s}{s}}\ell L^{-1/4}+ L^{-1/4}\mathrm{i}\log L +(\mp1+ \sqrt{3}\mathrm{i})w,\quad w\in\realR_{\ge 0}, 
\end{equation}
or
\begin{equation}
\zeta=\frac{X}{2(1-s)}\mp\sqrt{L}\pm\frac{1}{2}\sqrt{\frac{1-s}{s}}\ell L^{-1/4}- L^{-1/4}\mathrm{i}\log L +(\mp1- \sqrt{3}\mathrm{i})w,\quad w\in\realR_{\ge 0}, 
\end{equation}
we have
\begin{equation}
\label{eq:aux_10}
\begin{split}
&\pm\mathrm{Re}\left(g_2(\zeta)+\left(\frac{X^3}{24(1-s)^2}-\frac{XL_2}{2(1-s)}\right)\right)\\
&=\mathrm{Re}\left(\mp\frac{1-s}{3}\left(\zeta-\frac{X}{2(1-s)}\right)^3\pm L_2\left(\zeta-\frac{X}{2(1-s)}\right)\right)\\
&=-\frac{2}{3}(1-s)L^{3/2}+\sqrt{s(1-s)}\ell L^{3/4} -\frac{1}{4}s\ell^2-\frac{1}{4}sx^2 -(1-s)(\log L)^2+\bigO(L^{-3/4}(\log L)^3)\\
&\quad -\frac{8(1-s)}{3}w^3 -(2(1-s)L^{1/2}-\sqrt{s(1-s)}\ell L^{-1/4} +2\sqrt{3}(1-s)L^{-1/4}\log L)w^2\\
&\quad
-\left(2\sqrt{3}(1-s)L^{1/4}\log L+(1-s)L^{-1/2}(\log L)^2 -\sqrt{3s(1-s)}\ell L^{-1/2}\log L \right.\\
&\qquad\left.+\frac{x^2s}{4}L^{-1/2}-\frac{s}{4}\ell^2 L^{-1/2}\right)w.
\end{split}
\end{equation}
We remark that the only differences between these two formulas~\eqref{eq:aux_09},~\eqref{eq:aux_10} and two previous ones~\eqref{eq:aux_07},~\eqref{eq:aux_08} are the switch of $s$ and $1-s$, and the change of the sign before $\ell$.

Inserting these formulas~\eqref{eq:aux_07},~\eqref{eq:aux_08},~\eqref{eq:aux_09} and~\eqref{eq:aux_10} in~\eqref{eq:aux_06}, we know that the main contribution of the integral comes from the contours $\Gamma_{*}^{(L)}$. Also note that 
\begin{equation}
\frac{\left(\xi^{(1)}-\eta^{(2)}\right)\left(\eta^{(1)}-\xi^{(2)}\right)}{\left(\xi^{(1)}-\eta^{(1)}\right)\left(\xi^{(2)}-\eta^{(2)}\right)}\approx -1
\end{equation}
when all the variables are close to $-\sqrt{L}$ (for $\xi^{(i)}$) and $\sqrt{L}$ (for $\eta^{(i)}$). Here we use the notation $\approx$ for the leading order asymptotics. We obtain, 
after changing the integration variables,
\begin{equation}
\begin{split}
&\oint_0\ddbar{\rz}{(1-\rz)^2}\mathrm{T}_{1,1}(\rz;L_1,L_2,X;s)\\
&\approx \frac{e^{-\frac{4}{3}L^{3/2}-\frac12\ell^2-\frac12x^2}}{16\pi^4L}\\
&\quad\cdot\iiiint_{[-\log L,\log L]^4}e^{-s(u^{(1)})^2-s(v^{(1)})^2-(1-s)(u^{(2)})^2-(1-s)(v^{(2)})^2}\mathrm{d}u^{(1)}\mathrm{d}v^{(1)}\mathrm{d}u^{(2)}\mathrm{d}v^{(2)} \\
&\approx \frac{e^{-\frac{4}{3}L^{3/2}-\frac12\ell^2-\frac12x^2}}{16\pi^2s(1-s)L}.
\end{split}
\end{equation}
This proves Lemma~\ref{lm:01}.

\subsection{Proof of Lemma~\ref{lm:02}}

We fix three positive constants $c_1,c_2,c_3$ satisfying $c_1<c_2<c_3<2c_1$. We change the integration variables and deform the contours as follows
\begin{equation}
\begin{split}
&\xi^{(1)}_{i_1}=-\sqrt{L}-c_3L^{-1/4} +L^{-1/4}e^{\pm2\pi\mathrm{i}/3}u^{(1)}_{i_1}, \\
&\qquad \text{for}\quad \xi^{(1)}_{i_1}\in \Gamma_{\LL,\inn}=-\sqrt{L}-c_3L^{-1/4} +L^{-1/4}e^{\pm2\pi\mathrm{i}/3}\realR_{\ge 0},\\
&\xi^{(2)}_{i_2}=-\sqrt{L}-c_2L^{-1/4} +L^{-1/4}e^{\pm2\pi\mathrm{i}/3}u^{(2)}_{i_2},\\
& \qquad \text{for}\quad \xi^{(2)}_{i_2}\in\Gamma_{\LL}=-\sqrt{L}-c_2L^{-1/4} +L^{-1/4}e^{\pm2\pi\mathrm{i}/3}\realR_{\ge 0},\\
&\xi^{(1)}_{i_1}=-\sqrt{L}-c_1L^{-1/4} +L^{-1/4}e^{\pm2\pi\mathrm{i}/3}u^{(1)}_{i_1}, \\
&\qquad \text{for}\quad \xi^{(1)}_{i_1}\in\Gamma_{\LL,\out}=-\sqrt{L}-c_1L^{-1/4} +L^{-1/4}e^{\pm2\pi\mathrm{i}/3}\realR_{\ge 0},\\
&\eta^{(1)}_{i_1}=\sqrt{L}+c_3L^{-1/4} +L^{-1/4}e^{\pm\pi\mathrm{i}/3}v^{(1)}_{i_1},\\
& \qquad \text{for}\quad \eta^{(1)}_{i_1}\in\Gamma_{\RR,\inn}=\sqrt{L}+c_3L^{-1/4} +L^{-1/4}e^{\pm\pi\mathrm{i}/3}\realR_{\ge 0},\\
&\eta^{(2)}_{i_2}=\sqrt{L}+c_2L^{-1/4} +L^{-1/4}e^{\pm\pi\mathrm{i}/3}v^{(2)}_{i_2},\\
&\qquad \text{for}\quad \eta^{(2)}_{i_2}\in\Gamma_{\RR}=\sqrt{L}+c_2L^{-1/4} +L^{-1/4}e^{\pm\pi\mathrm{i}/3}\realR_{\ge 0},\\
&\eta^{(1)}_{i_1}=\sqrt{L}+c_1L^{-1/4} +L^{-1/4}e^{\pm\pi\mathrm{i}/3}v^{(1)}_{i_1}, \\
&\qquad \text{for}\quad \eta^{(1)}_{i_1}\in\Gamma_{\RR,\out}=\sqrt{L}+c_1L^{-1/4} +L^{-1/4}e^{\pm\pi\mathrm{i}/3}\realR_{\ge 0}.
\end{split}
\end{equation}
Note that these contours are nested in the order as in the definition.

Recall the definition of $g_1$ and $g_2$ in~\eqref{eq:def_g}. We show that with the above change of variables and the deformed contours,
\begin{equation}
\label{eq:estimate}
\begin{split}
&\mathrm{Re}\left(g_1(\xi^{(1)}_{i_1})\right)\le -\frac{2(1-\epsilon/2)}{3}sL^{3/2} +m(u^{(1)}_{i_1};s),\quad \text{for all}\quad\xi^{(1)}_{i_1}\in \Gamma_{\LL,\inn}\cup\Gamma_{\LL,\out},\\
&\mathrm{Re}\left(g_2(\xi^{(2)}_{i_2})\right)\le -\frac{2(1-\epsilon/2)}{3}(1-s)L^{3/2}+m(u^{(2)}_{i_2};1-s),\quad \text{for all}\quad\xi^{(2)}_{i_2}\in \Gamma_{\LL},\\	
&\mathrm{Re}\left(-g_1(\eta^{(1)}_{i_1})\right)\le -\frac{2(1-\epsilon/2)}{3}sL^{3/2}+m(v^{(1)}_{i_1};s),\quad \text{for all}\quad\eta^{(1)}_{i_1}\in \Gamma_{\RR,\inn}\cup\Gamma_{\RR,\out},\\
&\mathrm{Re}\left(-g_2(\eta^{(2)}_{i_2})\right)\le -\frac{2(1-\epsilon/2)}{3}(1-s)L^{3/2}+m(v^{(2)}_{i_2};1-s),\quad \text{for all}\quad\eta^{(2)}_{i_2}\in \Gamma_{\RR}
\end{split}
\end{equation}
for any fixed $\epsilon>0$ and sufficiently large $L$, where
\begin{equation}
m(w;s):=-sw^2.
\end{equation}

The proof of these inequalities are similar so we only prove the first one.

We drop the scripts and write $\xi=-\sqrt{L} -cL^{-1/4} + w(-1\pm\sqrt{3}\mathrm{i})L^{-1/4}$, where $w\in\realR_{\ge 0}$, and $c\in\{c_1,c_3\}$. Recall the scaling~\eqref{eq:scaling}. We have
\begin{equation}
\begin{split}
&\mathrm{Re}\left(g_1(\xi)\right)\\
&=-\frac{s}{3}\mathrm{Re}\left(\left(-\sqrt{L}-(c+w)L^{-1/4}+\sqrt{3}wL^{-1/4}\mathrm{i}\right)^3\right)\\
&\quad-\frac{X}{2}\mathrm{Re}\left(\left(-\sqrt{L}-(c+w)L^{-1/4}+\sqrt{3}wL^{-1/4}\mathrm{i}\right)^2\right)\\
&\quad +\left(L_1-\frac{X^2}{4s}\right)\mathrm{Re}\left(-\sqrt{L}-(c+w)L^{-1/4}+\sqrt{3}wL^{-1/4}\mathrm{i}\right)\\
&=-\frac{2}{3}sL^{3/2}+\bigO(L^{3/4})+(-2sw^2+2scw-\sqrt{s(1-s)}(x+\ell)w)\\
&\quad
+\left(-\frac{8}{3}sw^3-2scw^2+x\sqrt{s(1-s)}w^2 +(sc-\sqrt{s(1-s)})cw\right)L^{-3/4}\\
&=-\frac{2}{3}sL^{3/2}+m(w;s)+\bigO(L^{3/4})+(-sw^2+2scw-\sqrt{s(1-s)}(x+\ell)w)\\
&\quad
+\left(-\frac{8}{3}sw^3-2scw^2+x\sqrt{s(1-s)}w^2 +(sc-x\sqrt{s(1-s)})cw\right)L^{-3/4}
\end{split}
\end{equation}
where the term $\bigO(L^{3/4})$ is independent of $w$. Note that when $w\ll L^{3/8}$, the last two $w$ terms in the above expression are bounded by $L^{3/4}$. When $w\ge \bigO(L^{3/8})$, the last two $w$ terms behave like $-sw^2-\frac{8}{3}sw^3L^{-3/4}$ which is negative. Thus we have $\mathrm{Re}\left(g_1(\xi)\right)\le -\frac{2(1-\epsilon/2)}{3}sL^{3/2}+m(w;s)$ for sufficiently large $L$ and any $\xi\in\Gamma_{\LL,\inn}\cup\Gamma_{\LL,\out}$.

We need two more estimates before we prove Lemma~\ref{lm:02}. They were obtained in \cite{Liu21} so we do not provide the proof. The first one is (see (3.12) of \cite{Liu21})
\begin{equation}
\label{eq:aux_11}
\begin{split}
&
\left|\prod_{\ell=1}^2
\frac{\left( \Delta(\boldsymbol{\xi}^{(\ell)}) \right)^2
	\left( \Delta(\boldsymbol{\eta}^{(\ell)}) \right)^2
}
{\left( \Delta(\boldsymbol{\xi}^{(\ell)}; \boldsymbol{\eta}^{(\ell)}) \right)^2}
\cdot 
\frac{\Delta(\boldsymbol{\xi}^{(1)};\boldsymbol{\eta}^{(2)})
	\Delta(\boldsymbol{\eta}^{(1)};\boldsymbol{\xi}^{(2)})
}
{\Delta(\boldsymbol{\xi}^{(1)};\boldsymbol{\xi}^{(2)})
	\Delta(\boldsymbol{\eta}^{(1)};\boldsymbol{\eta}^{(2)})
}\right|\\
&\le k_1^{k_1/2}k_2^{k_2/2}(k_1+k_2)^{(k_1+k_2)/2}\prod_{i_1=1}^{k_1}\frac{1}{\dg(\xi_{i_1}^{(1)})}\frac{1}{\dg(\eta_{i_1}^{(1)})}\prod_{i_2=1}^{k_2}\frac{1}{\dg(\xi_{i_2}^{(2)})}\frac{1}{\dg(\eta_{i_2}^{(2)})}
\end{split}
\end{equation}
where  $\dg(\zeta)$ is the distance between $\zeta\in\Gamma_*$ and the contours $\Gamma_{\LL,\inn}\cup\Gamma_\LL\cup\Gamma_{\LL,\out}\Gamma_{\RR,\inn}\cup\Gamma_\RR\cup\Gamma_{\RR,\out}\setminus\Gamma_*$. In our choice of contours, it is easy to see that $\dg(\zeta)= \frac{\sqrt{3}}{2}cL^{-1/4}$ for $c\in\{c_3-c_2,c_2-c_1\}$. Hence we have
\begin{equation}
\label{eq:aux_13}
\begin{split}
&\left|\prod_{\ell=1}^2
\frac{\left( \Delta(\boldsymbol{\xi}^{(\ell)}) \right)^2
	\left( \Delta(\boldsymbol{\eta}^{(\ell)}) \right)^2
}
{\left( \Delta(\boldsymbol{\xi}^{(\ell)}; \boldsymbol{\eta}^{(\ell)}) \right)^2}
\cdot 
\frac{\Delta(\boldsymbol{\xi}^{(1)};\boldsymbol{\eta}^{(2)})
	\Delta(\boldsymbol{\eta}^{(1)};\boldsymbol{\xi}^{(2)})
}
{\Delta(\boldsymbol{\xi}^{(1)};\boldsymbol{\xi}^{(2)})
	\Delta(\boldsymbol{\eta}^{(1)};\boldsymbol{\eta}^{(2)})
}\right|\\
&\le k_1^{k_1/2}k_2^{k_2/2}(k_1+k_2)^{(k_1+k_2)/2} \left(\frac{\sqrt{3}}{2}\min\{c_3-c_2,c_2-c_1\}\right)^{-2k_1-2k_2}L^{(k_1+k_2)/2}
\end{split}
\end{equation}

The second estimate is, see (3.13) of \cite{Liu21},
\begin{equation}
\label{eq:aux_12}
\left|\mathrm{H}(\boldsymbol{\xi}^{(1)},\boldsymbol{\eta}^{(1)};
\boldsymbol{\xi}^{(2)},\boldsymbol{\eta}^{(2)})\right|\le \prod_{i_1=1}^{k_1}h(|\xi_{i_1}^{(1)}|)h(|\eta_{i_1}^{(1)}|)\prod_{i_2=1}^{k_2}h(|\xi_{i_2}^{(2)}|)h(|\eta_{i_2}^{(2)}|),
\end{equation}
where $h(y):=(1+y+y^2+y^3)^4$.

Inserting the estimates~\eqref{eq:estimate},~\eqref{eq:aux_13} and~\eqref{eq:aux_12} in~\eqref{eq:def_limit_T}, we obtain
\begin{equation}
\label{eq:aux_14}
\begin{split}
&\left|\mathrm{T}_{k_1,k_2}(\rz;L_1,L_2,X;s)\right|\\
&\le k_1^{k_1/2}k_2^{k_2/2}(k_1+k_2)^{(k_1+k_2)/2} \left(\frac{\sqrt{3}}{2}\min\{c_3-c_2,c_2-c_1\}\right)^{-2k_1-2k_2}L^{(k_1+k_2)/2} \\
&\quad\cdot  e^{-\frac{4}{3}(1-\epsilon/2)(sk_1+(1-s)k_2)L^{3/2}} \cdot    
\left|1-\mathrm{z}\right|^{k_2}
\left|1-\frac{1}{\mathrm{z}}\right|^{k_1}\\
&\quad
\prod_{i_1=1}^{k_1}     
\left(  \frac{1}{|1-\mathrm{z}|} \int_{\Gamma_{\LL,\inn}} 	    
\frac{|\mathrm{d}\xi_{i_1}^{(1)}|}{2\pi}  
+\frac{|\mathrm{z}|}{|1-\mathrm{z}|} \int_{\Gamma_{\LL,\out}} \frac{|\mathrm{d}\xi_{i_1}^{(1)}|}{2\pi}  
\right)e^{m(u_{i_1}^{(1)};s)}h(|\xi_{i_1}^{(1)}|)\\
&\quad
\prod_{i_1=1}^{k_1}   \left(  \frac{1}{|1-\mathrm{z}|} \int_{\Gamma_{\RR,\inn}} 	    
\frac{|\mathrm{d}\eta_{i_1}^{(1)}|}{2\pi}  
+\frac{|\mathrm{z}|}{|1-\mathrm{z}|} \int_{\Gamma_{\RR,\out}} \frac{|\mathrm{d}\eta_{i_1}^{(1)}|}{2\pi}  
\right)e^{m(v_{i_1}^{(1)};s)}h(|\eta_{i_1}^{(1)}|)
\\
&\quad \cdot 
\prod_{i_2=1}^{k_2}  
\int_{\Gamma_\LL}\frac{|\mathrm{d}\xi_{i_2}^{(2)}|}{2\pi}e^{m(u_{i_2}^{(2)};1-s)}h(|\xi_{i_2}^{(2)}|)\int_{\Gamma_\RR} \frac{|\mathrm{d}\eta_{i_2}^{(2)}|}{2\pi} 
e^{m(v_{i_2}^{(2)};1-s)}h(|\eta_{i_2}^{(2)}|).
\end{split}	
\end{equation}
On the other hand, note that the function $e^m$ decays super-exponentially along the integration contours. We obtain
\begin{equation}
\int_{\Gamma_{\LL,\inn}\cup\Gamma_{\LL,\out}} 	 e^{m(u_{i_1}^{(1)};s)}h(|\xi_{i_1}^{(1)}|)\frac{|\mathrm{d}\xi_{i_1}^{(1)}|}{2\pi} \le CL^{6}\cdot L^{-1/4}<CL^6
\end{equation} 
where $C$ is a positive constant. Here the term $L^6$ comes from the function $h$ near the point $u_{i_1}^{(1)}=0$, and $L^{-1/4}$ from the change of variables. We have similar estimates for all other integrals on the right hand side of~\eqref{eq:aux_14}. Thus we obtain
\begin{equation}
\begin{split}
&\left|\mathrm{T}_{k_1,k_2}(\rz;L_1,L_2,X;s)\right|\\
&\le k_1^{k_1/2}k_2^{k_2/2}(k_1+k_2)^{(k_1+k_2)/2} C^{k_1+k_2} L^{13(k_1+k_2)} e^{-\frac{4}{3}(1-\epsilon/2)(sk_1+(1-s)k_2)L^{3/2}}
\end{split}
\end{equation}
for a different constant $C$ if $|\rz|=r\in(0,1)$ is fixed. Note that $L^{13(k_1+k_2)}\ll e^{\frac{2}{3}\epsilon(sk_1+(1-s)k_2) L^{3/2}}$ as $L$ becomes large. Lemma~\ref{lm:02} follows immediately.

\def\cydot{\leavevmode\raise.4ex\hbox{.}}

\end{document}